\input amstex
\font\tensans=cmss10

\def\lj{\text{\tensans[}}
\def\rj{\text{\tensans]}}
\def\bbR{{\bold R}}
\def\bbZ{{\bold Z}}
\def\lh{L}

\def\z{{w}}
\def\hs{\hskip.7pt}
\def\hh{\hskip.4pt}
\def\nh{\hskip-.7pt}

\def\Om{\varOmega\hh}
\def\ric{\text{\rm Ric}}
\def\scal{\text{\rm Scal}}
%% begin section numbers
\def\id{0}
\def\as{1}
\def\pt{2}
\def\fc{3}
%% end section numbers
%% begin reference numbers
\def\besse{1}
\def\bourguignon{2}
\def\dm{3}
\def\fik{4}
\def\friedan{5}
\def\hamilton{6}
\def\ik{7}
\def\ivey{8}
\def\makowsky{9}
\def\perelman{10}
%% end reference numbers
\def\fmccl{1}
\def\betti{2}
\def\berge{3}
\def\pos{1}
\def\lrc{2}
\def\bai{3}
\def\lkc{4}
\documentstyle{amsppt}
\magnification=1200
\NoBlackBoxes
\topmatter
\title A Myers\hs-type theorem and compact Ricci solitons\endtitle
\rightheadtext{A Myers-type theorem and compact Ricci solitons}
\author Andrzej Derdzinski\endauthor
\address Department of Mathematics, Ohio State University, Columbus, OH 43210, 
USA\endaddress
\email andrzej\@math.ohio-state.edu\endemail
\abstract Let the Ricci curvature of a compact Riemannian manifold be 
greater, at every point, than the Lie derivative of the metric with respect to 
some fixed smooth vector field. It is shown that the fundamental group then 
has only finitely many conjugacy classes. This applies, in particular, to all 
compact shrinking Ricci solitons.\endabstract
\endtopmatter
\voffset=-35pt
\hoffset=15pt
\document 
\head\S\id. Introduction\endhead
Myers's classical theorem \cite{\besse, Theorem~6.51} implies that any compact 
Riemannian manifold with positive Ricci curvature has a finite fundamental 
group. Its standard proof uses a diameter estimate. In this note we observe 
that an easy upper bound \thetag{\lkc} on the length spectrum leads to a 
seemingly weaker but similar assertion, valid in a more general case:
\proclaim{Theorem~\fmccl}If\/ $\,(M,g)\,$ is a compact Riemannian manifold 
such that
$$\Cal L_\z g\,\,+\,\,\ric\,\,>\,\,0\tag\pos$$
for some $\,C^\infty$ vector field\/ $\,\z$, then\/ $\,\pi_1M\,$ has only 
finitely many conjugacy classes.
\endproclaim
Here $\,\Cal L_\z$ and $\,\hs\ric\hs\,$ are the Lie derivative and the 
Ricci tensor, while positivity in \thetag{\pos} means positive-definiteness at 
every point. As $\,H_1(M,\bbZ)=\pi_1M/\lj\hs\pi_1M,\pi_1M\hs\rj$, the 
following corollary is immediate:
\proclaim{Corollary~\betti}For any compact Riemannian manifold\/ $\,(M,g)\,$ 
satisfying the assumptions of Theorem~\fmccl, the homology group 
$\,H_1(M,\bbZ)\,$ is finite and\/ $\,b_1(M)=0$.
\endproclaim
It is not known if the conclusion of Theorem~\fmccl\ is actually weaker 
than finiteness of $\,\pi_1M$. Whether {\it a finitely presented group with 
only finitely many conjugacy classes must itself\hskip4ptbe finite\/} is an 
open question, raised by Makowsky \cite{\makowsky} in 1974. Another part of 
Myers's theorem, stating that any complete Riemannian manifold with 
$\,\hs\ric\,\ge c\hs g>0\,$ for a constant $\,c\,$ is necessarily 
compact, fails when $\,\hs\ric\hs\,$ is replaced by 
$\,\Cal L_\z g+\ric\hs$. Namely, Feldman, Ilmanen and Knopf 
\cite{\fik, Theorem~1.5} provide examples of noncompact 
complete K\"ahler manifolds $\,(M,g)\,$ with real holomorphic vector fields 
$\,\z\,$ satisfying the relation
$$\Cal L_\z g\,\,+\,\,\ric\,\,=\,\,c\hs g\hskip13pt\text{\rm for\ a\ 
constant}\hskip7ptc\hh,\tag\lrc$$
and in those particular examples $\,c>0$, so that 
$\,\Cal L_\z g+\ric\,\ge c\hs g>0$.

A Riemannian manifold $\,(M,g)\,$ admitting a $\,C^\infty$ vector field 
$\,\z\,$ with \thetag{\lrc} is called a {\it Ricci soliton\/} \cite{\fik} -- 
\cite{\hamilton}, \cite{\ivey}, \cite{\perelman}; one then also refers to 
$\,g\,$ as a {\it quasi-Einstein metric}. Such $\,(M,g)\,$ is said to be a 
{\it shrinking soliton\/} if $\,c>0\,$ in \thetag{\lrc}. 

The assumptions (and hence conclusions) of Theorem~\fmccl\ and 
Corollary~\betti\ clearly hold when $\,(M,g)\,$ is a compact shrinking 
Ricci soliton, or a small perturbation thereof. For compact shrinking Ricci 
solitons $\,(M,g)$, some special instances of this fact are known: finiteness 
of $\,H_1(M,\bbZ)\,$ is a trivial consequence of Theorem~1 of \cite{\ik} 
(see \S\fc); $\,M\,$ is simply connected if $\,g\,$ is also a K\"ahler 
metric (\S\fc); and in Proposition~2.2.5 on p.~396 of \cite{\friedan} it is 
stated that $\,b_1(M)=0\,$ under the additional assumption of a scalar 
curvature bound $\,\hs\scal\,\ge(\dim M-2)\hh c\hs$, for $\,c\,$ 
as in \thetag{\lrc}.

It is known that every non-Einstein compact Ricci soliton must be a shrinking 
soliton (cf.\ the end of \S\fc). This leads to a further conclusion:
\proclaim{Corollary~\berge}The Euler characteristic of every non-flat compact 
four-dimensional Ricci soliton is positive.
\endproclaim
In other words, Berger's inequality $\,\chi(M)>0$, for the Euler 
characteristic of a non-flat compact four-dimensional Einstein manifold 
\cite{\besse, Theorem~6.32}, remains true for Ricci solitons $\,(M,g)$. In 
fact, $\,\chi(M)=2+b_2(M)\,$ in the (orientable) non-Einstein case, since the 
soliton is shrinking and so, as stated above, $\,b_1(M)=0$.

The author wishes to thank Tadeusz Januszkiewicz, Gideon Maschler, Ronald 
Solomon and Fangyang Zheng for helpful comments.

\head\S\as. A simple estimate\endhead
Given $\,a,b\in\bbR\,$ with $\,a<b$, let a constant-speed geodesic 
$\,[a,b\hh]\ni t\mapsto x(t)\,$ in an orientable Riemannian manifold 
$\,(M,g)\,$ be smoothly closed, so that $\,x(b)=x(a)\,$ and 
$\,\dot x(b)=\dot x(a)$, and have the minimum length compared to all nearby 
smoothly closed $\,C^\infty$ curves $\,[a,b\hh]\to M$. Then
$$(b-a)\,\int_a^{\hs b}\ric\hs(\dot x,\dot x)\,dt\hskip8pt
\le\hskip8pt2\hs\lj\hs\text{\rm dist}\,(I,\varPsi)\hs\rj^2\hskip8pt
\le\hskip8ptk\pi^2\nh,\tag\bai$$
where $\,k\,$ is the largest even integer with $\,k+1\le\dim M$, while 
$\,\varPsi\in\,\text{\rm SO}\hs(m)$, for $\,m=\dim M-1$, is the {\it holonomy 
matrix\/} of our smoothly closed geodesic, characterized by the matrix-product 
relation $le\,\lj\hs e_1(a)\,\ldots\,e_m(a)\rj
=\hs\lj\hs e_1(b)\,\ldots\,e_m(b)\rj\hs\varPsi\,$ for some (or any) system 
$\,[a,b\hh]\ni t\mapsto e_j(t)\in T_{x(t})M\,$ of $\,m\,$ orthonormal vector 
fields parallel along the geodesic and normal to it. In addition, 
$\,I\in\,\text{\rm SO}\hs(m)\,$ stands for the identity matrix, and 
$\,\hs\text{\rm dist}\hs\,$ denotes the geodesic distance function in 
$\,\hs\text{\rm SO}\hs(m)\,$ corresponding to its submanifold metric induced 
by the inner product $\,Q\,$ in the ambient vector space 
$\,\hs\frak{gl}\hs(m,\bbR)\,$ of all $\,m\times m\,$ real matrices, given 
by $\,2\hs Q(\varPsi,\varPhi)=\,\text{\rm tr}\hskip3pt\varPsi\varPhi^*\nh$.

We now verify \thetag{\bai}. In any Riemannian manifold $\,(M,g)$, if 
$\,t\mapsto w(t)\,$ is a $\,C^\infty$ unit vector field normal to a fixed 
geodesic $\,[a,b\hh]\ni t\mapsto x(t)\,$ having the properties listed in the 
lines preceding \thetag{\bai}, and $\,w(b)=w(a)$, then 
$\,(R(\dot x,w)\dot x,w)\le(\dot w,\dot w)$, where 
$\,(w,w\hh')=\int_a^b\langle w,w\hh'\rangle\,dt\,$ stands for the $\,L^2$ 
inner product of vector fields $\,w,w\hh'$ tangent to $\,M\,$ along the 
geodesic, $\,R\,$ is the curvature tensor, and $\,\dot w=\nabla_{\!\dot x}w$. 
This is a well-known consequence of the length-minimizing property of the 
geodesic; cf.\ \cite{\ik, formula (3.1)}. Let us now select the 
$\,m\,$ fields $\,e_j$ as in the lines following \thetag{\bai}. For any fixed 
$\,C^\infty$ curve 
$\,[a,b\hh]\ni t\mapsto\varPhi(t)\in\,\text{\rm SO}\hs(m)\,$ joining $\,I\,$ 
to the holonomy matrix $\,\varPsi$, we may apply the inequality 
$\,(R(\dot x,w)\dot x,w)\le(\dot w,\dot w)\,$ to each of the 
$\,m\,$ fields $\,w=w_j$ given by the matrix-product formula 
$\,\lj\hs w_1(t)\,\ldots\,w_m(t)\rj
=\hs\lj\hs e_1(t)\,\ldots\,e_m(t)\rj\hs\varPhi(t)$. Summing the resulting 
inequalities over $\,j=1,\dots,m$, we get 
$\,\int_a^b\ric\hs(\dot x,\dot x)\,dt\le2\int_a^b|\dot\varPhi|^2\,dt$, with 
$\,|\hskip3pt|\,$ corresponding to the inner product $\,Q$. Since the curve 
$\,t\mapsto\varPhi(t)\,$ in $\,\hs\text{\rm SO}\hs(m)\,$ was arbitrary, we may 
choose it to be a constant-speed minimizing geodesic joining $\,I\,$ to 
$\,\varPsi$. The last inequality then yields the first relation in \thetag{\bai}.

To obtain the remaining inequality in \thetag{\bai}, it suffices to write any 
$\,\varPsi\in\,\text{\rm SO}\hs(m)\,$ as 
$\,\varPhi(\theta_1,\dots,\theta_p)$, which acts as the identity on 
$\,\lj\varPi_1\oplus\ldots\oplus\varPi_p\rj^\perp$, and as a rotation by the 
angle $\,\theta_l\in[-\pi,\pi]\,$ on each $\,\varPi_l$, for some fixed set 
of mutually orthogonal planes $\,\varPi_l$ in $\,\bbR^m$, $\,l=1,\dots,p$. The 
length of the curve 
$\,[\hs0,1\hs]\ni t\mapsto\varPhi(t\theta_1,\dots,t\theta_p)$, joining $\,I\,$ 
to $\,\varPsi$, then is 
$\,\lj\theta_1^2\nh+\dots+\theta_p^2\rj^{1/2}\le\sqrt{k/2\,}\hs\pi$.
\remark{Remark}The upper bounds in \thetag{\bai} reflect the following 
easily-verified facts, which we do not use: the metric in 
$\,\hs\text{\rm SO}\hs(m)\,$ corresponding to $\,\hs\text{\rm dist}\hs\,$ is 
bi-invariant, the diameter of $\,\hs\text{\rm SO}\hs(m)\,$ relative to 
$\,\hs\text{\rm dist}\hs\,$ equals $\,\sqrt{k/2\,}\hs\pi$, for $\,k\,$ as in 
\thetag{\bai}, and the curve joining $\,I\,$ to $\,\varPsi\,$ in 
$\,\hs\text{\rm SO}\hs(m)$, defined above, is a minimizing geodesic.
\endremark

\head\S\pt. Proof of Theorem~\fmccl\endhead
Give a compact, orientable Riemannian manifold $\,(M,g)\,$ with a $\,C^\infty$ 
vector field $\,\z\,$ satisfying \thetag{\pos}, let us choose a constant 
$\,c>0\,$ with $\,\Cal L_\z g+\,\ric\,\ge c\hs g$. We then have the 
following upper bound on the length spectrum:
$$\lh\,\,\le\,\,\sqrt{k/c\,}\,\pi\hskip8pt\text{\rm for\ the\ largest\ even\ 
integer}\hskip6ptk\,\le\,\dim M-1\hs,\tag\lkc$$
$\lh\,$ being the length of any smoothly closed constant-speed geodesic 
$\,[a,b\hh]\ni t\mapsto x(t)\,$ in $\,(M,g)\,$ which represents a local 
minimum of the length functional in its free homotopy class. Namely, 
\thetag{\bai} with $\,\Cal L_\z g+\,\ric\,\ge c\hs g\,$ gives 
$\,c\hs\lh^{\nh2}=(b-a)\hs c\int_a^bg(\dot x,\dot x)\,dt
\,\le\,(b-a)\int_a^b\ric\hs(\dot x,\dot x)\,dt
\,\le\,k\pi^2$, and \thetag{\lkc} follows. (The Lie-de\-riv\-a\-tive term in 
\thetag{\lrc} does not contribute to the integral, as 
$\,(\Cal L_\z g)(\dot x,\dot x)=2\hs d\hs\lj\hs g(\z,\dot x)\rj/dt$.)

The uniform bound \thetag{\lkc} implies in turn that there are only finitely 
many free homotopy classes of closed curves in $\,M$. In fact, an infinite 
sequence of smoothly closed geodesics $\,\gamma_j$ with uniformly bounded 
lengths $\,\lh_j$ cannot represent infinitely many distinct free homotopy 
classes, as one sees choosing a point $\,x_j$ on each $\,\gamma_j$ with a unit 
vector $\,u_j$ tangent to $\,\gamma_j$ at $\,x_j$, and then selecting a 
convergent subsequence of the sequence $\,(x_j,u_j,\lh_j)$.

Since the free homotopy classes are in a bijective correspondence with the 
conjugacy classes in the fundamental group of $\,M$, this completes the proof.

\head\S\fc. Comments\endhead
In this section we elaborate on some comments made in the introduction.

Ricci solitons on a compact manifold $\,M\,$ are precisely the fixed points of 
the Ricci flow $\,dg/dt=-\hs2\,\ric\hs\,$ projected, from the space of 
metrics, onto its quotient under diffeomorphisms and scalings 
\cite{\hamilton}. For shrinking solitons those scalings cause the metric to 
shrink to zero in finite time. Theorem~1 of \cite{\ik} states that, under the 
Ricci flow with any initial compact Riemannian manifold $\,(M,g)$, the lengths 
of curves representing a fixed element of infinite order in $\,H_1(M,\bbZ)\,$ 
remain bounded away from zero. As a consequence, if $\,(M,g)\,$ is a compact 
shrinking Ricci soliton, $\,H_1(M,\bbZ)\,$ has no element of infinite order, 
and so it is finite.

Secondly, compact shrinking Ricci solitons $\,(M,g)\,$ in which $\,g\,$ is a 
K\"ahler metric are known to be simply connected. In fact, $\,\z\,$ in 
\thetag{\lrc} then is holomorphic (cf.\ \cite{\dm}), and so the real 
cohomology classes of the Ricci and K\"ahler forms of $\,g\,$ are related by 
$\,[\rho]=c\hs[\Om]$, as \thetag{\lrc} gives $\,\rho=c\hh\Om-d\hh\xi\,$ for 
the $\,1$-form $\,\xi=\imath_\z\Om$. Thus, $\,c_1(M)>0$, and so 
$\,\pi_1M=\{0\}\,$ by a result of Kobayashi \cite{\besse, Theorem~11.26}.

Finally, $\,c>0\,$ for every non-Einstein compact Ricci soliton $\,(M,g)$. 
In fact, by Ivey's Proposition~1 in \cite{\ivey}, the scalar curvature 
$\,\hs\scal\hs\,$ of $\,g\,$ then must be positive, while taking the 
$\,g$-trace of \thetag{\lrc} we see that $\,c\,$ is the average value of 
$\,\hs\scal\hs\,$ on $\,M$. (That $\,c>0\,$ and $\,\hs\scal\hs\ge0\,$ is 
also stated in Propositions~2.2.2 -- 2.2.4 on p.~396 of \cite{\friedan}, 
while $\,\hs\scal\hs\,$ is nonconstant, as shown by Bourguignon 
\cite{\bourguignon, Proposition~3.11}.)

\enddocument